\documentclass{amsproc}

\usepackage{a4wide}
\usepackage{amsmath}
\usepackage{enumerate}

\usepackage{amsfonts, amssymb,amsthm,amsgen}

\def\Ric{\mathop{{\rm Ric}}\nolimits}
\def\tRic{\mathop{\widetilde{\rm Ric}}\nolimits}
\def\ker{\mathop{{\rm ker}}\nolimits}

\def\Real{\mathbb{R}}
\def\Co{\mathbb{C}}
\def\g{\mathfrak{g}}
\def\h{\mathfrak{h}}

\def\so{\mathfrak{so}}
\def\spin{\mathfrak{spin}}
\def\simil{\mathfrak{sim}}

\def\su{\mathfrak{su}}
\def\u{\mathfrak{u}}
\def\z{\mathfrak{z}}
\def\sp{\mathfrak{sp}}
\def\f{\mathfrak{f}}
\def\zr{\ltimes}
\def\pr{\mathop\text{\rm pr}\nolimits}
\def\tr{\mathop\text{\rm tr}\nolimits}
\def\rk{\mathop\text{\rm rk}\nolimits}

\def\p{\partial}
\def\P{\mathcal{P}}
\def\R{\mathcal{R}}
\def\be{\begin{equation}}
\def\ee{\end{equation}}
\def\E{{\mathcal E}}

\newtheorem{theorem}{Theorem}
\newtheorem{lem}{Lemma}
\newtheorem{prop}{Proposition}

\begin{document}

\title[How to find Lorentzian holonomy algebra]
{How to find the holonomy algebra of a Lorentzian manifold}

\author{Anton S. Galaev}
%\thanks{University of Hradec Kr\'alov\'e, Faculty of Science, Jana Koziny
%1237, 500~03 Hradec Kr\'alov\'e,  Czech Republic\\
%E-mail: anton.galaev(at)uhk.cz}

\maketitle

\begin{abstract}
Manifolds with exceptional holonomy play an important role in
string theory, supergravity and M-theory. It is explained how one
can find the holonomy algebra of an arbitrary Riemannian or
Lorentzian manifold. Using the de~Rham and Wu decompositions, this
problem is reduced to the case of locally indecomposable
manifolds. In the case of locally indecomposable Riemannian
manifolds, it is known that the holonomy algebra can be found from
the analysis of special geometric structures on the manifold. If
the holonomy algebra $\mathfrak{g}\subset\mathfrak{so}(1,n-1)$ of
a locally indecomposable  Lorentzian manifold $(M,g)$ of dimension
$n$ is different from $\mathfrak{so}(1,n-1)$, then it is contained
in the similitude algebra $\mathfrak{sim}(n-2)$. There are 4 types
of such holonomy algebras. Criterion how to find the type of
$\mathfrak{g}$ are given,
 and special geometric structures corresponding to each type are described.
To each $\mathfrak{g}$ there is a canonically associated
subalgebra $\mathfrak{h}\subset\mathfrak{so}(n-2)$. An algorithm
how to find $\mathfrak{h}$ is provided.

{\bf Keywords:} Lorentzian manifold, holonomy group, holonomy
algebra,  de~Rham-Wu decomposition

{\bf MSC 2010  codes:} 53C29, 53C25, 53C50, 81T30

\end{abstract}

\section{Introduction}
Riemannian manifolds with special holonomy play an important role
in theoretical physics, in particular, in
 string theory compactifications, supergravity and M-theory, see
\cite{Cecotti,CGLP,Gubser,Gu-Sp} and references therein. Very
important are explicit constructions of complete or compact
Riemannian manifolds with special holonomy, since they give
ansatzs to physical theories. The main reason is that these
manifolds are Ricci-flat and admit parallel spinors.

 To find the holonomy algebra of an indecomposable Riemannian manifolds,
one may use the fact that Riemannian  manifolds with different
holonomy algebras have different geometric properties, in
particular, they admit different parallel forms, see
\cite{Besse,Joyce07} and Section \ref{secRiemhol} below. Another
approach using the curvature tensor can be found in \cite{McIn}.
In Section \ref{secdeR} we explain how to find the holonomy
algebra of an arbitrary Riemannian manifold using our algorithm
for the de~Rham decomposition from \cite{Gal14}.

Holonomy algebras of 4-dimensional Lorentzian manifolds and their
relation to General Relativity were studied e.g. in
\cite{KG1,H-L,JJ}. The classification of the holonomy algebras of
Lorentzian manifolds of arbitrary dimension was achieved in
\cite{BB-I,Leistner,Galmetr}. This classification is quite
different from the one in the Riemannian case. The only
irreducible holonomy algebra is $\so(1,n+1)$, ${\rm dim} M=n+2\geq
3$. Other holonomy algebras of locally indecomposable manifolds
are contained in the similitude Lie algebra
$\simil(n)=(\Real\oplus\so(n))\zr\Real^n$, and this is the most
interesting case.

Recently an attention to this classification was taken by
theoretical physicists. In \cite{BCH1} it is noted that compering
with Riemannian manifolds with special  holonomy, the Lorentzian
ones have many interesting and unusual properties in the context
of string theories. In \cite{G-P} the Einstein equation on the
manifolds with $\simil(n)$-holonomy was considered, and the
partial solution were interpreted as multi-centered black holes.
In \cite{CGHP} the universality of $\simil(n)$-holonomy metrics as
the solutions to the Einstein and supergravity equations were
studied. In \cite{Cecotti} the F-theory on $\simil(n)$-holonomy
manifolds was constructed. The historical review \cite{Gibbons09}
provides a bridge from Riemannian manifolds with special holonomy
to the Lorentzian ones and it gives a motivation for study
Lorentzian metrics with special holonomy. See also
\cite{Coh-Gl,Grover}. Parallel spinors on Lorentzian manifolds and
their relation to holonomy were studied in
\cite{Bryant,FF00,Leistner}; more general equations on spinors are
considered in the review \cite{Baum08}. In \cite{B-M,Bazajkin}
constructions of globally hyperbolic Lorentzian manifolds with
special holonomy are given. Examples of Lorentzian manifolds with
different holonomy and some global geometrical properties will
provide ansatzs to physical theories. In particular, this leads to
the problem to find a way to compute the holonomy algebra of an
arbitrary Lorentzian manifold.

In this paper we describe geometrical properties of Lorentzian
manifolds with different types of holonomy algebras and we give a
complete answer to the natural question: {\it How to find the
holonomy algebra of an arbitrary Lorentzian manifold?}

To a holonomy algebra $\g\subset\simil(n)$ one associates  its
$\so(n)$-projection $\h$, which must be the holonomy algebra of a
Riemannian manifold \cite{Leistner} and it is called the
orthogonal part of $\g$. We provide an algorithm how to find $\h$
in Section \ref{secorthpart}.

Next,  there are 4 types of holonomy algebras $\g\subset\simil(n)$
corresponding to locally indecomposable manifolds. The algebras of
type 1 and 2 have simple structure and they are of the form
$(\Real\oplus\h)\zr\Real^n$ and $\h\zr\Real^n$, respectively. The
algebras of types 3 and 4 are more exotic, and they can be
obtained from the first two by some twistings. In Section
\ref{sectype} we provide criteria that allow to find the type of
the holonomy algebra. Similar criteria are given in \cite{Boubel}.
Our criteria are more concrete: on one hand, we show how the type
of the holonomy algebra can be found using the local coordinates
and it becomes computable, on the other hand, we describe the
geometric structures corresponding to each of the types.

In Section \ref{sectype} we explain also how to use our algorithm
from \cite{Gal14}  for obtaining the Wu decomposition of
Lorentzian manifolds in order to reduce the computation of the
holonomy algebra of an arbitrary Lorentzian manifold to the case
of locally indecomposable one.

Thus we provide the complete algorithm that allows to find the
holonomy algebra of an arbitrary Riemannian or Lorentzian
manifold. This algorithm can be computerized, e.g. as a part of
the package Differential Geometry for Maple \cite{Anderson}, since
it requires computations of certain parallel tensors (e.g.
parallel symmetric bilinear forms and certain  differential forms
that can be found as the solutions to some systems of partial
differential equations), some computations in linear algebra, and
computations in local coordinates.

Note that we consider only the holonomy algebra, i.e. we study the
connected component of the holonomy group. The full holonomy group
of a not simply connected manifold can be bigger. There are
several results about this group in the Riemannian case
\cite{Besse,McIn1,McIn2}, and a recent paper \cite{BLL} for the
case of Lorentzian manifolds.

\section{Holonomy group; holonomy algebra}

The theory of holonomy groups of pseudo-Riemannian manifolds can
be found e.g. in \cite{Besse,Joyce07}. Let $(M,g)$ be a Riemannian
or Lorentzian manifold. {\it The holonomy group} $G_x$ of $(M,g)$
at a point $x\in M$ is the  Lie group that consists of the
orthogonal (resp. Lorentzian) transformations of the tangent space
given by the parallel transports along all piecewise smooth loops
at the point $x$.  The corresponding Lie algebra  is called {\it
the holonomy algebra} and it determines the holonomy group if the
manifold is simply connected. {\it The Ambrose-Singer Theorems}
states that the holonomy algebra is spanned by the endomorphisms
 $\tau^{-1}_\gamma\circ R_y(\tau_\gamma X,\tau_\gamma
Y)\circ \tau_\gamma$ of $T_xM$, where  $\gamma$ is a piecewise
smooth curve starting at the point $x$ with an end-point $y \in
M$, and $X,Y\in T_xM$.

 The {\it fundamental principle for holonomy groups}
states that there exists a one-to-one correspondents between
parallel tensor fields $T$ ($\nabla T=0$) on $M$ and tensors $T_0$
of the same type at $x$ preserved by the tensor extension of the
representation of the holonomy group.

Since we are interested in holonomy algebras, {\it in what follows
we will assume that the manifold $M$ is simply connected}. In
general case one can pass to the universal covering $(\tilde
M,\tilde g)$.

In general it is impossible to find the holonomy group using the
definition, and it is impossible to find the holonomy algebra
using the Ambrose-Singer Theorem, since then one should consider
parallel transports along all piecewise smooth loops at a point or
parallel transports along all piecewise smooth curves starting at
a point. Below we will show how to compute the holonomy algebra of
any Riemannian or Lorentzian manifold. For that we will use the
classification of the holonomy algebras for these manifolds and
the geometric properties of manifolds with each possible holonomy
algebra.

\section{The de~Rham and Wu decompositions}\label{secdeR}

The de~Rham and Wu decompositions allow to decompose a Riemannian
or Lorentzian manifold into a local product of locally
indecomposable manifolds. In the case of Riemannian manifolds, the
local indecomposability is equivalent to the irreducibility of the
holonomy group. In the case of  Lorentzian manifolds, the local
indecomposability is equivalent to the weak irreducibility of the
holonomy group (a subgroup of the Lorentzian group is called
weakly irreducible if it does not preserve any non-degenerate
proper subspace of the Minkowski space).

In the original proofs of the de~Rham and Wu theorems it is
supposed that the holonomy group is know. This makes unclear, how
to find the holonomy group of an arbitrary manifold. By this
reason in \cite{Gal14} we give algorithms for finding the de~Rham
and Wu decompositions for Riemannian and Lorentzian manifolds. The
algorithms use the analysis of parallel symmetric bilinear forms
on the manifold $(M,g)$, i.e., we do not require the knowledge of
the holonomy group. More precisely, we find parallel distributions
$E_0,...,E_r$ on $M$ such that the induced connection on $E_0$ is
flat, and the holonomy group of $(M,g)$ is the product of the
holonomy groups of the induced connections on the distributions
$E_1,...,E_r$, that are (weakly-) irreducible.

\section{Riemannian holonomy algebras}\label{secRiemhol}

The results that we review in this section are the major
achievements of the holonomy theory, they can be found e.g. in
\cite{Besse,Joyce07}. We will need them in Section
\ref{secorthpart}.

The holonomy algebra of a locally symmetric Riemannian space
$(M,g)$ at a point $x$ coincides with $\{R_x(X,Y)|X,Y\in T_xM\}$.
Locally $(M,g)$ is isometric to a symmetric space $H/G$, where $H$
is the group of transvections of that space; the holonomy group of
that space coincides with the isotropy representation of the
stabilizer $H$ of a point. The list of indecomposable Riemannian
spaces can be found e.g. in \cite{Besse}.

Here we list irreducible holonomy algebras of not locally
symmetric Riemannian manifolds and we give the description of the
corresponding geometries (on a simply connected manifold)
including the Einstein condition and parallel forms (we do not
include trivial parallel forms i.e. constant function and the
volume form on an orientable manifold):

\begin{itemize}
\item $\so(n)$: generic Riemannian manifolds, no parallel forms;
\item $\u(m)\subset\so(2m)$: K\"ahlerian manifolds, parallel K\"ahlerian 2-form and its powers,
not Ricci-flat;
\item $\su(m)\subset\so(2m)$: special K\"ahlerian manifolds or Calabi-Yau manifolds, parallel K\"ahlerian 2-form, its
powers, parallel complex volume form and its conjugate,
Ricci-flat;
\item $\sp(k)\subset\so(4k)$: hyper-K\"ahlerian manifolds, 3 independent parallel K\"ahlerian 2-forms and forms obtained
from their combinations,
 Ricci-flat;
\item $\sp(k)\oplus\sp(1)\subset\so(4k)$: quaternionic-K\"ahlerian manifolds, parallel 4-from and its powers,
Einstein and not
Ricci-flat;
\item $\spin(7)\subset\so(8)$: Ricci-flat, parallel 4-form;
\item $G_2\subset\so(7)$: Ricci-flat, a parallel 3-form and its dual.
\end{itemize}

Irreducible Riemannian holonomy algebras $\g\subset\so(n)$ that
appear as the holonomy algebras of symmetric Riemannian spaces and
are different from $\so(n)$, $\u(m)$,  $\sp(k)\oplus\sp(1)$ are
called {\it symmetric Berger algebras}.

Now we may easily find the holonomy algebra $\g$ of any Riemannian
manifold $(M,g)$. Recall that passing to the universal covering,
we may assume that $M$ is simply connected. First suppose that
$(M,g)$ is locally indecomposable, i.e. its holonomy algebra is
irreducible. If $\nabla R=0$, then the manifold is locally
symmetric and its holonomy algebra at a point $x\in M$ coincides
with $\{R_x(X,Y)|X,Y\in T_xM\}$. If $\nabla R\neq 0$, then there
are only 7 possibilities for $\g$. According to the list of the
holonomy algebras and to the geometric properties of the
corresponding Riemannian manifolds, to find $\g$ it is enough to
compute the Ricci tensor of $(M,g)$ and to find parallel
2,3,4-forms on $(M,g)$ (all that can be done using e.g. Maple); of
course, one should also analyze the dimension $n$ of $(M,g)$, e.g.
if $n=9$, then $\g=\so(9)$; if $n=7$, then $\g=G_2$ if and only if
there exists a parallel 3-form  and $\g=\so(7)$ otherwise. Another
approach that uses the computation of the curvature tensor is
proposed in \cite{McIn}.

Consider now an arbitrary simply connected Riemannian manifold
$(M,g)$.
 Results from \cite{Gal14} allow us to to find the distributions
$E_\alpha$, $1\leq \alpha\leq r$, defining the Wu decomposition of
$(M,g)$. The holonomy algebra of $(M,g)$ is the direct sum of the
holonomy algebras of the induced connections on the distributions
$E_\alpha$, $1\leq \alpha\leq r$, considered as  vector bundles
over $M$, see Section \ref{secdeR}. The holonomy algebra of each
distribution is irreducible and it can be found in the same way as
the holonomy algebra of a locally indecomposable Riemannian
manifold above. Let $\nabla$ and $R$ be the connection and the
curvature of the manifold $(M,g)$. If $\nabla R|_{E_\alpha\times
E_\alpha\times E_\alpha}=0$, then the holonomy algebra of
$E_\alpha$ coincides with $\{R_x(X,Y)|X,Y\in E_{\alpha x}\}$;
otherwise it is one of the holonomy algebras from the above list
and it can be found analyzing $\Ric|_{E_\alpha\times E_\alpha}$
and parallel sections of the bundle $\Lambda^k E_\alpha$,
$k=2,3,4$.

\section{Classification of the Lorentzian holonomy algebra }\label{secLorHol}

Here we review results from \cite{BB-I,Leistner,Galmetr,ESI}. Let
$(M,g)$ be a simply connected Lorentzian manifold of dimension
$n+2$, $n\geq 0$. Fix a point $x\in M$. The tangent space
$(T_xM,g_x)$ can be identified with the Minkowski space
$(\Real^{1,n+1},g_x)$.  Then the holonomy algebra  $(M,g)$ at the
point $x$ is identified with a subalgebra $\g\subset\so(1,n+1)$.
From the above it follows that we may assume that
$\g\subset\so(1,n+1)$ is weakly irreducible. If
$\g\subset\so(1,n+1)$ is irreducible, then $\g=\so(1,n+1)$.
Suppose that $\g\subset\so(1,n+1)$ is not  irreducible, the
$\g$ preserves an isotropic line in $\Real^{1,n+1}$.

 We fix a basis $p,X_1,...,X_n,q$ of $\Real^{1,n+1}$ such that $p$ and $q$ are light-like vectors,
 $g(p,q)=1$ and  the subspace $E\subset\Real^{1,n+1}$
 spanned by $X_1,...,X_n$ is an Euclidean subspace orthogonal to $p$ and $q$. We obtain the decomposition
\begin{equation}\label{decTxM} T_xM=\Real^{1,n+1}=\Real p\oplus E\oplus \Real q.\end{equation}
Denote  by $\simil(n)$ the subalgebra of $\so(1,n+1)$ that preserves the isotropic line $\Real p$. The Lie algebra
$\simil(n)$ can be identified with the following matrix algebra:
 \begin{equation}\label{matsim}\simil(n)=\left\{\left. \left (\begin{array}{ccc} a
&(GX)^t & 0\\ 0 & A &-X \\ 0 & 0 & -a \\
\end{array}\right)\right|\, a\in \Real,\, A \in \so(n),\, X\in \Real^n \right\},\end{equation}
where $G$ is the Gram matrix of the restriction of $g$ to $E$ with
respect to the basis $X_1,\dots,X_n$. The above matrix can be
identified with the triple $(a,A,X)$. We get the decomposition
$$\simil(n)=(\Real\oplus\so(n))\zr\Real^n,$$ which means that $\Real\oplus\so(n)\subset\simil(n)$ is a subalgebra and
$\Real^n\subset\simil(n)$ is an ideal, and the Lie brackets of
$\Real\oplus\so(n)$ with $\Real^n$ are given by the standard
representation of $\Real\oplus\so(n)$ in $\Real^n$. The Lie
algebra $\simil(n)$ is isomorphic to the Lie algebra of the Lie
group of similarity transformations of $\Real^n$.  We assume that
$\g\subset\simil(n)$. We identify $\Real^n$ and $E$.

Let $\h\subset\so(n)$ be a subalgebra. Recall that $\h$ is a
compact Lie algebra and we have the decomposition
$\h=\h'\oplus\z(\h)$, where  $\h'$ is the commutant of $\h$, and
$\z(\h)$ is the center of $\h$. If $\h\subset\so(n)$ is
irreducible, then $\z(\h)\neq 0$ implies
$\h\subset\u(\frac{n}{2})$; in this case
$\h'\subset\su(\frac{n}{2})$ and $\z(\h)=\Real J$, where $J$ is
the complex structure.

The next theorem gives the classification of weakly irreducible
not irreducible holonomy algebras of Lorentzian manifolds.

\begin{theorem}
A subalgebra $\g\subset \simil(n)$ is the weakly irreducible
holonomy algebra of a Lorentzian manifold if and only if it is
conjugated to one of the following subalgebras:
\begin{description}
\item[type 1.] $\g^{1,\h}=(\Real\oplus\h)\zr\Real^n$,

\item[type 2.] $\g^{2,\h}=\h\zr\Real^n$,

\item[type 3.] $\g^{3,\h,\varphi}=\{(\varphi(A),A,0)|A\in\h\}\zr\Real^n,$

\item[type 4.] $\g^{4,\h,m,\psi}=\{(0,A,X+\psi(A))|A\in\h,X\in\Real^m\}$,  \end{description}

where $\h\subset\so(n)$ is the holonomy algebra of a Riemannian
manifold; for $\g^{3,\h,\varphi}$ it holds $\z(\h)\neq\{0\}$, and
$\varphi :\h\to\Real$ is a non-zero linear map with
$\varphi|_{\h'}=0$; for $\g^{4,\h,m,\psi}$ it holds
  $2\leq m<n$ is an integer,
$\h\subset\so(m)$,  $\dim\z(\h)\geq n-m$, a decomposition
$\Real^n=\Real^m\oplus \Real^{n-m}$ is fixed, and $\psi:\h\to
\Real^{n-m}$ is a surjective linear map with $\psi|_{\h'}=0$.

\end{theorem}

The subalgebra $\h\subset\so(n)$ associated to a weakly
irreducible Lorentzian holonomy algebra $\g\subset \simil(n)$  is
called {\it the orthogonal part} of $\g$. For $\h\subset\so(n)$
there exist the decompositions
\begin{equation}\label{LM0A}E=E_0\oplus E_1\oplus\cdots\oplus E_r,\quad
\h=\{0\}\oplus\h_1\oplus\cdots\oplus\h_r\end{equation} such that
$\h$ annihilates $E_0$, $\h_\alpha(E_\beta)=0$ for $\alpha\neq
\beta$, and $\h_\alpha\subset\so(E_\alpha)$ is an irreducible
subalgebra for $1\leq \alpha\leq r$. Let $n_\alpha={\rm dim} E_\alpha$.

Let us give a more precise descriptions for algebras of type 3.
Let $K\subset\{1,...,r\}$ be the set of indices such that
$\varphi|_{\h_\alpha}\neq 0$. If $\alpha\in K$, then $\h_\alpha
\subset\u(\frac{n_\alpha }{2})$; in this case
$\h'_\alpha\subset\su(\frac{n_\alpha}{2})$ and
$\z(\h_\alpha)=\Real J_\alpha $, where $J_\alpha $ is  the complex
the structure on $E_\alpha$.  Let $c_\alpha=\varphi(J_\alpha)$.
Then
\begin{equation}\label{type3A}
\g^{3,\h,\varphi}=\big(\oplus_{\alpha\in K}
\Real(c_\alpha+J_\alpha)\oplus\h'_\alpha\,\oplus\,\,\oplus_{\alpha\not\in
K}\h_\alpha \big)  \zr\Real^n,\end{equation} where
$c_\alpha+J_\alpha$ denotes $(c_\alpha,J_\alpha,0).$

Similarly we may write \begin{equation}\label{type4}
\g^{4,\h,m,\psi}=\big(\oplus_{\alpha\in K}
\Real(J_\alpha+\psi(J_\alpha))\oplus\h'_\alpha\,\oplus\,\,\oplus_{\alpha\not\in
K}\h_\alpha \big)  \zr\Real^m,\end{equation} where $K$ is defined
in the same way as for $\g$ of type 3.

\section{Walker and adapted coordinates}\label{secWalker}

Let $(M,g)$ be a Lorentzian manifold  with the holonomy algebra
$\g\subset\simil(n)$. Then $(M,g)$ admits a parallel distribution
of isotropic lines $\ell$.  According to \cite{Walker}, locally there
exist so called Walker coordinates $v,x^1,...,x^n,u$ such that the
metric $g$ has the form
\begin{equation}\label{Walker} g=2dvdu+h+2Adu+H (d
u)^2,\end{equation} where $h=h_{ij}(x^1,...,x^n,u)d x^id x^j$ is
an $u$-dependent family of Riemannian metrics,  $A=A_i(x^1,
\ldots, x^n,u)d x^i$ is an $u$-dependent family of one-forms, and
$H=H(v,x^1,...,x^n,u)$ is a local function on $M$. An important
example of such spaces are pp-waves that are given by $h=\sum
(dx^i)^2$, $A=0$, $\partial_vH=0$. Equivalently, pp-waves are
Walker manifolds with the holonomy algebras contained in
$\Real^n\subset\simil(n)$.

Consider the local frame
$$p=\p_v,\quad X_i=\p_i-A_i\p_v,\quad q=\p_u-\frac{1}{2}H\p_v.$$
Let $E$ be the distribution generated by the vector fields
$X_1$,...,$X_n$. Clearly, the vector fields $p$, $q$ are
light-like, $g(p,q)=1$, the restriction of $g$ to $E$ is positive
definite, and $E$ is orthogonal to $p$ and $q$. The vector field
$p$ defines the parallel distribution of isotropic lines $\ell$
and it is recurrent, i.e. $\nabla p=\mu\otimes p$, where
$\mu=\frac{1}{2}\p_vHdu$. Since the manifold is locally
indecomposable, any other recurrent vector field is proportional
to $p$. Next, $p$ is proportional to a parallel vector field if
and only if $d\theta =0$, which is equivalent to
$\p^2_vH=\p_i\p_vH=0$. In the last case the coordinates can be
chosen in such a way that $\p_v H=0$ and $\nabla p=\nabla \p_v=0$,
see e.g. \cite{ESI}.

  Boubel \cite{Boubel} proved  that there exist
Walker coordinates
$$v,\,x_0=(x_0^1,\ldots ,x_0^{n_0}),\ldots ,x_r=(x_r^1,...,x_r^{n_r}),\,u$$
adapted to the decomposition \eqref{LM0A}. This means that \be
\label{h=sum}h=h_0+h_1+\cdots+h_r,\quad h_0=\sum_{i=1}^{n_0}(d
x_0^i)^2,\quad h_\alpha=\sum_{i,j=1}^{n_\alpha}h_{\alpha ij} d
x_{\alpha}^id x_{\alpha}^j,\ee
$$A=\sum_{\alpha=0}^r A_\alpha,\quad A_0=0,\quad A_\alpha=\sum_{k=1}^{n_\alpha} A^\alpha_k d x^k_\alpha,$$
and   one has \be\label{halpha=} \frac{\p}{\p
{x^k_\beta}}h_{\alpha ij}=\frac{\p}{\p
{x^k_\beta}}A^\alpha_i=0,\quad  \text{ if }\beta\neq\alpha.\ee We
call these coordinates {\it adapted}.
 The coordinates can be chosen
 so that in addition $A=0$, see~\cite{GL10}.

\section{The curvature tensor}\label{secR}

Since the Ambrose-Singer Theorem provides the relation of the
holonomy algebra and the curvature tensor, we describe here the
curvature tensor of a Walker manifold $(M,g)$ with the holonomy
algebra $\g\subset\simil(n)$ at the point $x\in M$.  For that it
is convenient to consider the space $\R(\g)$ of algebraic
curvature tensors of type $\g$, i.e. the space of linear maps from
$\Lambda^2\Real^{1,n+1}$ to $\g$ satisfying the first Bianchi
identity. The curvature tensor of $(M,g)$, $R=R_x$ at the point
$x\in M$ belongs to the space $\R(\g)$. This space is found in
\cite{Gal1,onecomp}. Consider the decomposition \eqref{decTxM}.
For a subalgebra $\h\subset\so(n)$  consider the space
$$\P(\h)=\{P\in (\Real^n)^*\otimes
\h|g(P(X)Y,Z)+g(P(Y)Z,X)+g(P(Z)X,Y)=0\text{ for all } X,Y,Z\in
\Real^n\}.$$ Define the map $\tRic:\P(\h)\to\Real^n$,
$\tRic(P)=P^j_{ik}g^{ik}X_j$. It does not depend on the choice of
the basis $X_1,...,X_n$. The tensor $R\in\R(\g^{1,\h})$ is
uniquely given by elements $\lambda\in\Real,e\in
E,R_0\in\R(\h),P\in\P(\h),T\in\odot^2E$ in the following way:
\begin{align*}
R(p,q)=&(\lambda,0,e),\qquad R(X,Y)=(0,R_0(X,Y),P(Y)X-P(X)Y),\\
R(X,q)=&(g(e,X),P(X),T(X)),\qquad R(p,X)=0 \end{align*}
 for all $X,Y\in\Real^n$.
We write $R=R(\lambda,e,R_0,P,T)$. The Ricci tensor $\Ric(R)$ of
$R$ is given by  $\Ric(R)(U,V)={\rm tr}(Z\mapsto R(Z,U)V)$ and it
satisfies
\begin{align}\label{Ric1} \Ric(p,q)=&\lambda,\quad \Ric(X,Y)=\Ric(R_0)(X,Y),\\
\label{Ric2} \Ric(X,q)=&g(X,e-\tRic(P)),\quad \Ric(q,q)=-{\rm tr}
T.
\end{align}

Decomposition \eqref{LM0A} defines the decompositions
$P=P_1+\cdots+P_r$, $P_\alpha\in\P(\h_\alpha)$ and
$R_0=R_{01}+\cdots+P_{0r}$, $R_{0\alpha}\in\R(\h_\alpha)$.

For the above tensor $R$, the condition
$R\in\R(\g^{3,\h,\varphi})$ is equivalent to the following
conditions:
$$\lambda=0,\qquad g(e,X)=\varphi(P(X)),\quad X\in E,\qquad R_0\in\R(\ker\varphi).$$
The condition $R\in\R(\g^{4,\h,m,\psi})$ is equivalent to the
following conditions:
$$\lambda=0,\qquad e=0,\qquad \pr_{\Real^{n-m}}\circ T=\psi\circ P,  \qquad R_0\in\R(\ker\psi).$$

Note that  a weakly irreducible holonomy algebra
$\g\subset\simil(n)$
 defines canonically only the isotropic line $\Real p$.
 Let us take a real number $\mu\neq 0$, the vector $p'=\mu p$, and
any light-like vector $q'$ with $g(p',q')=1$. There exists a
unique vector $w\in E$ such that
$q'=\frac{1}{\mu}(-\frac{1}{2}g(w,w) p+w+q)$.
 The corresponding
$E'$ has the form $E'=\{-g(x,w)p+x|x\in E\}$. We will consider the
map $x\in E\mapsto x'=-g(x,w)p+x\in E'$. Using this, we  obtain
that $R=R(\tilde\lambda,\tilde e,\tilde R_0,\tilde P,\tilde T)$.
For example, it holds \be\label{changecomp}\tilde
\lambda=\lambda,\quad \tilde e=\frac{1}{\mu}(e-\lambda w)',\quad
\tilde P(x')=\frac{1}{\mu}(P(x)+R_0(x,w))',\quad \tilde
R_0(x',y')z'=(R_0(x,y)z)'.\ee This shows e.g. that if $\lambda=0$,
then the projection of the vector $e$ to $p^\perp/\Real p$ is
defined up to a non-zero real multiple.

Let $n=2m\geq 2$ and consider the space $\P(\u(m))$, note that
$\u(1)=\so(2)$. In \cite{onecomp} it is shown that the
$\u(m)$-module $\P(\u(m))$ admits the decomposition
$$\P(\u(m))=\P_0(\u(m))\oplus\P_1(\u(m))$$
into the direct sum of irreducible submodules. It holds
$\P_0(\u(m))=\{P\in\P(\u(m))|\tRic(P)=0\}$ and $\P_1(\u(m))\simeq
\Real^n$. The last isomorphism has the form
$$Z\in\Real^n\mapsto P,\quad P(X)=R^{CP^m}(X,Z),$$
where $R^{CP^m}$ is the curvature tensor at a point of the complex
projective space,
$$R^{CP^m}(X,Z)=\frac{1}{2}g(JX,Z)J+\frac{1}{4}(X\wedge Z+JX\wedge
JZ),$$
where $(X\wedge Z)Y=g(X,Y)Z-g(Z,Y)X$.

\begin{lem}\label{lemPZ}
For $P\in\P_1(\u(m))$ corresponding to $Z\in\Real^n$ it holds
$$\tRic(P)=\frac{m+1}{2}Z,\quad \pr_{\Real
J}P(X)=\frac{m+1}{2m}g(JX,Z)J.$$\end{lem}

{\it Proof.} Using the complex structure $J$, we identify the
space $\Real^{2m}$ with $\Co^m$. Let $\tilde g$ be the Hermitian
form on $\Co^m$ corresponding to $g$, i.e. $$\tilde
g(X,Y)=g(X,Y)+ig(X,JY).$$ Let $e_1,...,e_m$ be an orthogonal basis
of $\Co^m$. For the trace of any element $L\in \u(m)$ acting on
$\Co^m$ it holds
$$\tr_\Co L=\sum_{k=1}^m\tilde g(L e_k,e_k)=\sum_{k=1}^m(g(L e_k,e_k)+ig(L
e_k,Je_k))=i\sum_{k=1}^mg(L e_k,Je_k).$$ Recall that for $L\in
\su(m)$ it holds $\tr_\Co L=0$, and $\tr_\Co J=mi$. Note that
$$(X\wedge Z+JX\wedge
JZ)Y=\tilde g(Y,X)Z-\tilde g(Y,Z)X.$$ This implies that
\begin{multline*}\tr_\Co(X\wedge Z+JX\wedge
JZ)=\sum_{k=1}^m \tilde g((X\wedge Z+JX\wedge JZ)e_k,e_k)\\=
 \sum_{k=1}^m \tilde g(\tilde g(e_k,X)Z-\tilde g(e_k,Z)X,e_k)=\tilde
g\left(Z,\sum_{k=1}^m\tilde g(X,e_k),e_k\right)-\tilde
g\left(X,\sum_{k=1}^m\tilde g(Z,e_k),e_k\right)\\=\tilde
g(Z,X)-\tilde g(X,Z)=2i g(Z,JX).\end{multline*} We conclude that
$$\pr_{\Real J}P(X)=\pr_{\Real
J}R^{CP^m}(X,Z)=\frac{m+1}{2m}g(JX,Z)J.$$ In \cite{ESI} it is
shown that $$g(\tRic P,X)=-\sum_{k=1}^mg(P(JX)e_k,Je_k).$$ for all
$X\in E$. Consequently,
$$g(\tRic
P,X)=i\tr_\Co P(JX)=i\tr_\Co R^{CP^m}(JX,Z)=\frac{m+1}{2}g(X,Z),$$
i.e. $\tRic P=\frac{m+1}{2}Z.$ \hfill $\Box$

The above considerations easily imply the following

\begin{lem}\label{lemopr} Suppose that $\h\subset\u(\frac{n}{2})$. If
$\lambda=0$ and $\Ric(R_0)=0$, then the  projections of the
vectors $e$ and $Z$ to $p^\perp/\Real p$ are defined up to a
non-zero real multiple. \end{lem}

If we fix on $(M,g)$ Walker coordinates as in Section
\ref{secWalker}, then we get vector fields $p$, $q$, and  a
distribution $E$ over an open subset of $M$. Consequently, the
curvature tensor of $(M,g)$ over this subset is defined by some
tensor fields $\lambda,v,R_0,P,T$. It can be checked that
$R_0=R(h)$ is the curvature tensor of the Riemannian metric $h$,
and it holds
\begin{equation}\label{lamv}\lambda=-\frac{1}{2}\p^2_vH,\quad
e=-\frac{1}{2}\left(\p_i\p_vH-A_i\p^2_vH\right)h^{ij}X_j.\end{equation}

\section{Finding the orthogonal part of a Lorentzian holonomy algebra}\label{secorthpart}

Let $(M,g)$ be a simply connected Lorentzian manifold with the
holonomy algebra $\g\subset\simil(n)$. In this section we give an
algorithm how to find the orthogonal part $\h\subset\so(n)$ of
$\g$.

The subalgebra $\h\subset\so(n)$ coincides with the holonomy
algebra of the induced connection on the so-called screen bundle
$\E=\ell^\perp/\ell$ \cite{Leistner1}. If we choose a
decomposition \eqref{decTxM} over an open subset of $M$, then $\E$
restricted to this subset may be identified with the distribution
$E$. For the curvature tensor of the connection on $\E$ we get
$$R(p,\cdot)=0,\quad R(X,Y)=R_0(X,Y),\quad R(X,q)=P(X),\quad X,Y\in\Gamma(E).$$

Recall that  $\h\subset\so(n)$ is the holonomy algebra of a
Riemannian manifold. The decomposition
$$\E=\E_0\oplus\E_1\oplus\cdots\oplus\E_r$$ into the direct sum of a
flat subbundle $\E_0\subset\E$ and parallel subbundles
$\E_1,...,\E_r\subset\E$, corresponding to the decompositions
\eqref{LM0A}, can be obtained exactly in the same way as in
Section \ref{secdeR}.  Hence we may assume that the subalgebra
$\h\subset\so(n)$ is irreducible.

The subalgebra $\h\subset\so(n)$ can not be found exactly in the
same way, as the holonomy algebra of a Riemannian manifold in
Section \ref{secRiemhol}, since we can not distinguish  symmetric
Berger subalgebras from the very beginning. By this reason we use
a deeper analysis.

If we already know  that  $\h\subset\so(n)$ is a  symmetric Berger
subalgebra, then it can be found in the following way. Let $y\in
M$ be a point such that either $R_y(h)\neq 0$, or $P_y\neq 0$,
such point exists since $\h\neq 0$. Since $\h$ does not contain
any proper Berger subalgebra, and each of   the subsets $\{
R_y(h)(X,Y)|X,Y\in E_y\}$ and $\{P_y(X)|X\in E_y\}$ generates a
Berger subalgebra in $\h$, $\h$ is  generated either by $\{
R_y(h)(X,Y)|X,Y\in E_y\}$, or by $\{P_y(X)|X\in E_y\}$.

First we compute $\Ric(h)$ and $\tRic(P)$ (for the last object,
the formula \eqref{comptRic} given below can be used). From the
results of \cite{onecomp} it follows that if $\Ric(h)=0$ and
$\tRic(P)=0$, then the subalgebra $\h\subset\so(n)$ is one of
$\so(n)$, $\su\left(\frac{n}{2}\right)$,
$\sp\left(\frac{n}{4}\right)$, $\spin(7)$, $G_2$. In this case,
$\h\subset\so(n)$ can be found simply analyzing the parallel
sections of $\Lambda^k\E$ $(k=2,3,4)$, as in Section
\ref{secRiemhol}.

Now we may assume that $\Ric(h)\neq 0$ or $\tRic(P)\neq 0$, then
$\h\subset\so(n)$ is either one  of $\so(n)$,
$\u\left(\frac{n}{2}\right)$,
$\sp\left(\frac{n}{4}\right)\oplus\sp(1)$, or $\h\subset\so(n)$
is a symmetric Berger subalgebra.

If there exists a parallel section of $\Lambda^2\E$, then $\h$ is
contained in $\u\left(\frac{n}{2}\right)$. Next, the subspace of
the $\u\left(\frac{n}{2}\right)$-module
$$\odot^2\u\left(\frac{n}{2}\right)\simeq\odot^2\su\left(\frac{n}{2} \right)\oplus\su\left(\frac{n}{2}\right)\oplus\Real$$
annihilated by $\u\left(\frac{n}{2}\right)$ is clearly  of
dimension two. This subspace is spanned by the curvature tensor of
the complex projective space and by the subset $\Real$. On the
other hand, any symmetric Berger subalgebra
$\h\subset\u\left(\frac{n}{2}\right)$ annihilates in addition the
curvature tensor valued at a point of the  corresponding symmetric
space, which is an element of the space
$\odot^2\u\left(\frac{n}{2}\right)$. Consequently,  if the space
of parallel sections of $\odot^2\u(\E)$ equals 2, then
$\h=\u\left(\frac{n}{2}\right)$. Otherwise,
$\h\subset\u\left(\frac{n}{2}\right)$ is a symmetric Berger
subalgebra.

Now  $\h\subset\so(n)$ is either one  of $\so(n)$,
$\sp\left(\frac{n}{4}\right)\oplus\sp(1)$, or it is a symmetric
Berger subalgebra not contained in $\u\left(\frac{n}{2}\right)$.

In \cite{Kost56} it is shown that if an indecomposable  simply
connected Riemannian symmetric space admits a non-trivial parallel
4-form, then its holonomy algebra is not simple. This and the list
of indecomposable simply connected Riemannian symmetric space
\cite{Besse} show that such a space admits a parallel 4-form, then
it is either K\"ahlerian, or quaternionic-K\"ahlerian, or its
holonomy algebra is one of $\so(r)\oplus\so(s)\subset\so(rs)$
($r,s\neq 2$) and $\sp(r)\oplus\sp(s)\subset\sp(rs)$  ($r,s\neq
1$).

If there are no non-trivial parallel sections in  $\Lambda^4\E$,
then either $\h=\so(n)$, or $\h\subset\so(n)$ is a simple
symmetric Berger algebra, which is not contained neither  in
$\u\left(\frac{n}{2}\right)$, nor in
$\sp\left(\frac{n}{4}\right)\oplus\sp(1)$. The Lie algebra
$\so(n)$  annihilates exactly one 1-dimensional subspace of
$\odot^2\so(n)$ , while for symmetric Berger subalgebras
$\h\subset\so(n)$  this subspace is at least of dimension 2.
Hence, if the dimension of parallel sections of the bundle
$\wedge^2\so(\E)$ equals to 1, then $\h=\so(n)$. Otherwise
$\h\subset\so(n)$ is a symmetric Berger algebra.

Suppose that there is a non-trivial parallel sections  in
$\omega\in\Lambda^4\E$. Note that the stabilizer of the Kraines
4-form evaluated at a point of a quaternionic-K\"ahlerian manifold
of dimension $n$ coincides with
$\sp\left(\frac{n}{4}\right)\oplus\sp(1)$ \cite{Besse}. Let
$\f\subset\so(n)$ be the stabilizer of $\omega$ at some point.
Clearly, $\h\subset\f$. If  $\f\neq
\sp\left(\frac{n}{4}\right)\oplus\sp(1)$, then $\h\subset\so(n)$
is a symmetric Berger algebra and it is one of
$\so(r)\oplus\so(s)\subset\so(rs)$ ($r,s\neq 2$) and
$\sp(r)\oplus\sp(s)\subset\sp(rs)$  ($r,s\neq 1$). Otherwise,
$\h\subset\f=\sp\left(\frac{n}{4}\right)\oplus\sp(1)$. Again, the
$\sp\left(\frac{n}{4}\right)\oplus\sp(1)$-module
$$\odot^2\left(\sp\left(\frac{n}{4}\right)\oplus\sp(1)\right)\simeq\odot^2\sp\left(\frac{n}{4}\right)\oplus
\sp\left(\frac{n}{4}\right)\otimes\sp(1)\oplus\odot^2\sp(1)$$
annihilated by $\sp\left(\frac{n}{4}\right)\oplus\sp(1)$ is of
dimension two, while for a symmetric Berger subalgebra
 $\h\subset\sp\left(\frac{n}{4}\right)\oplus\sp(1)$  this dimension is bigger.
 It is  enough to find parallel sections of the bundle $\odot^2\left(\sp\left(\E\right)\oplus\sp(1)\right)$.

\section{Finding the type of a Lorentzian holonomy algebra}\label{sectype}

Let $(M,g)$ be a locally indecomposable simply connected
Lorentzian manifold with the holonomy algebra
$\g\subset\so(1,n+1)$. If $(M,g)$ admits a parallel distribution
of isotropic lines, or equivalently, locally it admits recurrent
light-like vector fields that are proportional on the
intersections of the domains of their definition, then
$\g\subset\simil(n)$; otherwise, $\g=\so(1,n+1)$.

 We consider the case  $\g\subset\simil(n)$.    The
  following statement has been already discussed in Section
  \ref{secWalker}.

\begin{prop}\label{Prop1234} Let  $(M,g)$ be a locally indecomposable
 Lorentzian manifold with the holonomy algebra $\g\subset\simil(n)$, then $\g$ is of
 type 2 or 4 if and only if for any Walker coordinate system it holds  $\partial^2_v H=\partial_v\partial_i H=0$, equivalently, there exists a Walker coordinate system in a neighborhood of each point such that $\partial_v H=0$. If $M$ is simply connected, then these conditions are equivalent to the existence of a parallel light-like vector field.  \end{prop}

Now we should be able to decide  between types 1 and 3 and between types 2 and 4. We will do that in the following two theorems.

Suppose that  $\g$ is of type 1 or 3. The following theorem allows
to find the type of $\g$.

\begin{theorem}\label{ThType3}  Let $(M,g)$ be a simply
connected locally indecomposable Lorentzian
 manifold with the holonomy algebra $\g\subset\simil(n)$. Suppose that $(M,g)$  does not admit any parallel light-like vector field. Then $\g$ is of type 3 if and only if the following conditions hold:

\begin{itemize}

\item[1.] for any Walker coordinate system it holds $\partial^2_vH=0$;

\item[2.] there is a non-empty subset $K\subset\{1,...,r\}$ of indexes $\alpha$  such that

\begin{itemize}

\item[2.a.]  if $\alpha\in K$, then  $\h_\alpha$ is contained in
$\u\left(\frac{n_\alpha}{2}\right)$, i.e. the bundle $\E_\alpha$
admits a parallel complex structure $J_\alpha$;

 \item[2.b.]  the Riemannian metric $h_\alpha$ is Ricci-flat for all  $\alpha\in K$;

\item[2.c.] if $\alpha\not\in K$, then $\partial_{x^i_\alpha}\partial_vH=0$
 for any adapted coordinate system;

\item[2.d.] for each $\alpha\in K$ there exists a non-zero constant
$c_\alpha\in\Real$ such that  for any adapted coordinate system it
holds

\be\label{ur2d}pr_{E_\alpha}e=-\frac{2c_\alpha}{n_\alpha}J_\alpha\tRic
P_\alpha,\ee where $$e=-R(p,q)q,\quad \tRic
P_\alpha=\sum_{i,j=1}^{n_\alpha}h_\alpha^{ij}R(\partial_{x^i_\alpha},q)\partial_{x^j_\alpha}$$
and $pr_{E_\alpha}e\neq 0$ for some adapted coordinate system.

Equation \eqref{ur2d} has the following coordinate form
\be\label{ur2dA}
\partial_{x^i_\alpha}\partial_vH=-\frac{2c_\alpha}{n_\alpha}
\left(\nabla^jF^\alpha_{lj}+\nabla^j\dot h_{\alpha
lj}-\partial_{x^l_\alpha}h_\alpha^{jk}\dot h_{\alpha
jk}\right)J^l_{\alpha i},\ee where the indexes $j,k,l$ run from
$1$ to $n_\alpha$, $\dot h=\partial_u h$,
$F^\alpha_{lj}=\partial_{x^l_\alpha}A^\alpha_k-\partial_{x^k_\alpha}A^\alpha_l$
(no sum over $\alpha$). It holds
$\partial_{x^i_\alpha}\partial_vH\neq 0$ for some adapted
coordinate system and some $i$, $1\leq i\leq n_\alpha$.
\end{itemize}
\end{itemize}

Otherwise, $\g$ is of type 1.
\end{theorem}

{\bf Proof of Theorem \ref{ThType3}.}  Suppose that the holonomy
algebra $\g$ of $(M,g)$ at a point $x\in M$ is of type 3, i.e.
$\g=\g^{3,\h,\varphi}$. Let $K\subset\{1,...,r\}$ be as in Section
\ref{secLorHol}. Condition 2.a follows from the definition of the
algebra $\g^{3,\h,\varphi}$.

Let $y\in M$ and let $\gamma$ be a piecewise smooth curve
beginning at $x$ and ending at $y$. Note that the holonomy algebra
at the point $y$ equals to
$\tau_\gamma^{-1}\circ\g\circ\tau_\gamma$ and it is isomorphic to
$\g=\g^{3,\h,\varphi}$. Clearly, for the curvature tensor at the
point $y$ it holds
$R_y\in\R(\tau_\gamma^{-1}\circ\g\circ\tau_\gamma)\simeq\R(\g)$.
We will identify $\tau_\gamma^{-1}\circ\g\circ\tau_\gamma$ and
$\g$.
 Fix a
coordinate system in a neighborhood of the point $y$. Then $R_y$
can be decomposed as in Section \ref{secR}. Condition 1 of the
theorem follows from \eqref{lamv} and the fact that $\lambda=0$
for any element from $\R(\g)$. Let $\alpha\in K$ and $X,Y\in
E_{\alpha y}$. We have
$R_0(X,Y)=R(h)(X,Y)=R(h_\alpha)(X,Y)\in\h_\alpha$. The fact that
$R_0=R(h)\in\R(\ker\varphi)$ implies condition 2.b.

Let $\alpha\not\in K$ and $X\in E_{\alpha y}$. Then
$R(X,q)=(g(e,X),P(X),T(X))$ and $P(X)\in\h_\alpha$. Since $\g$
contains $\h_\alpha$ and $\Real^n$, we obtain $(g(e,X),0,0)\in\g$.
Consequently, $g(e,X)=0$ for any $X\in E_{\alpha y}$. This and
\eqref{lamv} imply 2.c.

Let $\alpha\in K$ and $X\in E_{\alpha y}$. The projection of
$R(X,q)$ on $\Real\oplus\Real J_\alpha\subset\simil(n)$ must
belong to $\Real(c_\alpha+J_\alpha)$. On the other hand,
$R(X,q)=(g(e,X),P(X),T(X))$ and  $\pr_{\Real\oplus \Real
J_\alpha}=g(e,X)+\pr_{\Real J_\alpha}P_\alpha(X)$. From Lemma
\ref{lemPZ} it follows that $\pr_{\Real
J_\alpha}P_\alpha(X)=-\frac{2}{n_\alpha}g( X,J_\alpha\tRic
P_\alpha)J_\alpha$. We conclude that $\pr_{E_{\alpha
x}}e=-\frac{2c_\alpha}{n_\alpha}J_\alpha\tRic P_\alpha$.

Let us find the coordinate form of the last equality. For
simplicity we assume that $n_\alpha=n_1=n$, i.e.
$\h\subset\u(\frac{n}{2})$ is irreducible.
 Since $\partial^2_vH=0$, from \eqref{lamv} it follows that
$e=-\frac{1}{2}(\partial_i\partial_vH)h^{ij}X_j$. In Section
\ref{secR} we have seen that $\Ric(X,q)=g(e-\tRic P,X)$ for all
$X\in E$. In \cite{G-P} it is shown that
\be\label{Riciqpq}\Ric(\partial_i,q)=-\frac{1}{2}\left(\partial_i\partial_vH+\nabla^jF_{ij}+\nabla^j\dot
h_{ij}-\partial_ih^{jk}\dot h_{jk}\right),\quad
\Ric(p,q)=-\frac{1}{2}\partial^2_vH.\ee Recall that
$X_i=\partial_i-A_ip$. We obtain that \be\label{comptRic}\tRic
P=(\tRic P)_ih^{ij} X_j,\quad (\tRic
P)_i=\frac{1}{2}\left(\nabla^jF_{ij}+\nabla^j\dot
h_{ij}-\partial_{i}h^{jk}\dot h_{jk}\right).\ee The equation under
consideration takes the form
$$(\partial_i\partial_v H)h^{ij}X_j=\frac{2c_\alpha}{n_\alpha}J(\tRic
P)_ih^{ij}X_j.$$ Let $JX_j=J_i^jX_j$. Then
$$(\partial_i\partial_v H)h^{ij}X_j=\frac{2c_\alpha}{n_\alpha}(\tRic
P)_ih^{ij}J^l_jX_l.$$ This implies
$$\partial_i\partial_v H=\frac{2c_\alpha}{n_\alpha}(\tRic
P)_kh^{kj}J^l_jh_{li}.$$ Since $J$ is a K\"ahlerian structure, it
holds $J^l_jh_{li}=-h_{jl}J^l_i$. Now it is easy to obtain
\eqref{ur2dA}. Thus condition 2 is proved.

Conversely, suppose that for a Lorentzian manifold $(M,g)$ the
conditions 1 and 2 hold. We should prove that the holonomy algebra
of $(M,g)$ at a point $x\in M$ coincides with $\g^{3,\h,\varphi}$.
Let $\gamma$ be a piecewise smooth curve beginning at $x$; let
$y\in M$ be its end-point. Fix a decomposition \eqref{decTxM} of
$T_xM$.  It defines  the decomposition \be\label{decTyM}
T_yM=\Real \tau_\gamma p \oplus \tau_\gamma E_x\oplus \Real
\tau_\gamma q_x.\ee  Since we know that $\g$ is either of type 1
or of type 3, it contains the ideal $\Real^n$. Consider the tensor
$R_\gamma=\tau_\gamma^{-1}\circ
R_y(\tau_\gamma\cdot,\tau_\gamma\cdot)\circ\tau_\gamma\in\R(\g).$
As in Section \ref{secR}, it is defined by elements $\lambda$,
$e$, $P$, $R_0$ and $T$. Since we have the isomorphism
$\tau_\gamma:T_xM\to T_yM$ and consider the decomposition
\eqref{decTyM}, the tensor $R_y$ is defined by the above elements
mapped by the isomorphism $\tau_\gamma$ to the point $y$. Fix an
adapted coordinate system in a neighborhood of the point $y$. Let
$\lambda_y,...,T_y$ be the elements defining $R_\gamma$ and
corresponding to these coordinates. The condition 1 and the
results of Section \ref{secR} imply that $\lambda=\lambda_y=0$.
Suppose that $\alpha\not\in K$. Condition 2.c implies that
$\pr_{E_{y\alpha}}e_y=0$. From Lemma \ref{lemopr} it follows that
$\pr_{E_{x\alpha}}e=0$, i.e. $\pr_\Real R_\gamma(X,q_x)=0$ for any
$X\in E_{x\alpha}$. Condition 2.b implies that
$R_\gamma(X,Y)\in\h'_\alpha$ for all $X,Y\in E_{x\alpha}$. We have
only to consider the projection $\pr_{\Real\oplus\Real J_\alpha}
\{R_\gamma(X,q_x)| \alpha\in K,\, X\in E_{x\alpha}\}.$ Lemma
\ref{lemopr}, condition 2.d and the above proof show that this
projection coincides with $\Real(c_\alpha+J_\alpha)$. Thus
$\g=\g^{3,\h,\varphi}$. The theorem is true. \hfill $\Box$

Suppose that we have a local Walker metric that satisfies
conditions 2.a and 2.b. According to Section \ref{secorthpart},
these conditions depend only on $h$ and $A$. We may ask if the
function $H$  can be found in such a way that the holonomy of this
metric is of type 3. The condition 1 of the theorem  can be easily
satisfied and we are left with Equations \eqref{ur2dA}. The
integrability condition for this system of equations is of the
form $\partial_{x^s_\alpha}B_i=\partial_{x^i_\alpha}B_s$, where
$B_i$ is the right hand side of \eqref{ur2dA}. Thus a priory the
function $H$ can not be changed to make the holonomy of the metric
to be of type 3.

Next suppose that the holonomy algebra $\g$ is of type 2 or 4.
Suppose that $\E_0\neq 0$ (this is true if $\g$ is of type 4).
Since the connection on the vector bundle $\E_0$ is flat and $M$
is simply connected, there exist orthonormal parallel sections
$e_1,...,e_{n_0}$ spanning $\E_0$. {\it We will assume} that all
adapted coordinate systems are chosen in such a way that
$\partial_{x^s_0}=e_s$. Suppose that $\g=\g^{4,\h,m,\psi}$. Suppose that $e_1,...,e_{n-m}$ is a basis of
$\psi(\h)$. Then there
exists a matrix $(c_{s\alpha})$ such that
$\psi(J_\alpha)=\sum_{s=1}^{n-m}c_{s\alpha}e_s$, $\alpha\in K$.
Note that  $m=n-{\rm rk}(c_{s\alpha})$.

\begin{theorem}\label{ThType4}  Let $(M,g)$ be a simply
connected locally indecomposable Lorentzian
 manifold with the holonomy algebra $\g\subset\simil(n)$. Suppose that $(M,g)$   admits a
  parallel light-like vector field. Then $\g$ is of type 4 if and only if
there exists a number $m$, $2\leq m<n$ such that
  the following conditions hold:

\begin{itemize}

\item[1.]   rank of the subbundle $\E_0\subset \E$ is not smaller then $n-m$;
 for any adapted coordinate system it
holds $\partial_{x^s_0}\partial_{x^t_0}H=0$, $1\leq s,t\leq n-m$;

\item[2.] there is a non-empty subset $K\subset\{1,...,r\}$ of indexes $\alpha$  such that

\begin{itemize}

\item[2.a.]  if $\alpha\in K$, then  $\h_\alpha$ is contained in
$\u\left(\frac{n_\alpha}{2}\right)$, i.e. the bundle $\E_\alpha$
admits a parallel complex structure $J_\alpha$;

 \item[2.b.]  the Riemannian metric $h_\alpha$ is Ricci-flat for all  $\alpha\in K$;

\item[2.c.] if $\alpha\not\in K$, then $\partial_{x^i_\alpha}\partial_{x^s_0}H=0$
 for any adapted coordinate system and $1\leq s\leq n-m$;

\item[2.d.] for each $\alpha\in K$ there exist numbers
$(c_{s\alpha})_{s=1}^{n-m}\in\Real$ such that $m=n-{\rm
rk}(c_{s\alpha})$  and  for any adapted coordinate system it holds

\be\label{ur2d4}pr_{E_\alpha}T(e_s)=-\frac{2c_{s\alpha}}{n_\alpha}J_\alpha\tRic
P_\alpha,\ee where $$T(e_s)=-R(e_s,q)q,\quad \tRic
P_\alpha=\sum_{i,j=1}^{n_\alpha}h_\alpha^{ij}R(\partial_{x^i_\alpha},q)\partial_{x^j_\alpha}$$
and for each $s$, $1\leq s\leq n-m$, $pr_{E_\alpha}T(e_s)\neq 0$ for some adapted coordinate system.

Equation \eqref{ur2d4} has the following coordinate form
\be\label{ur2dA4}
\partial_{x^i_\alpha}\partial_{x^s_0}H=-\frac{2c_{s \alpha}}{n_\alpha}
\left(\nabla^jF^\alpha_{lj}+\nabla^j\dot h_{\alpha
lj}-\partial_{x^l_\alpha}h_\alpha^{jk}\dot h_{\alpha
jk}\right)J^l_{\alpha i},\ee where the indexes $j,k,l$ run from
$1$ to $n_\alpha$, $\dot h=\partial_u h$,
$F^\alpha_{lj}=\partial_{x^l_\alpha}A^\alpha_k-\partial_{x^k_\alpha}A^\alpha_l$
(no sum over $\alpha$). For each $s$, $1\leq s\leq n-m$, it holds
$\partial_{x^i_\alpha}\partial_{x^s_0}H\neq 0$ for some adapted
coordinate system and some $i$, $1\leq i\leq n_\alpha$.
\end{itemize}
\end{itemize}
Otherwise, $\g$ is of type 2.
\end{theorem}

{\bf Proof of Theorem \ref{ThType4}.} The proof is similar to the
proof of Theorem \ref{ThType3}. The description of an element
$R\in \R(\g^{4,\h,m,\psi})$ implies $\pr_{\Real^{n-m}}\circ
T|_{\Real^{n-m}}=0$, this gives the condition
$\partial_{x^s_0}\partial_{x^t_0}H=0$, $1\leq s,t\leq n-m$. Next,
if $\alpha\not\in K$, then $\pr_{\Real^{n-m}}\circ
T|_{E_\alpha}=0$, this gives condition 2.c. Let us consider
condition 2.d. As in the same way as in the proof of Theorem
\ref{ThType3} we get
$$\pr_{E_0}T(X)=-\frac{2}{n_\alpha}g(X,J_\alpha\tRic
P)\psi(J_\alpha),\quad \alpha\in K, X\in E_{\alpha x}.$$
Substituting $\psi(J_\alpha)=\sum_{t=1}^{n-m}c_{t\alpha}e_t$,
multiplying the obtained equality by $e_s$, and using the facts
that $T$ is symmetric and $e_1,...,e_{n_0}$ is an orthonormal
basis, we obtain \eqref{ur2d4}. \hfill $\Box$

Thus, in order to find the type of the holonomy algebra $\g$
knowing the orthogonal part of $\g$, it is enough to apply
Proposition \ref{Prop1234} and to check the conditions 1 and 2 of
one of Theorems \ref{ThType3} or \ref{ThType4}. For that it is
necessary to find the parallel complex structures $J_\alpha$ on
the bundles $\E_\alpha$, or the corresponding parallel 2-forms. To
compute the right hand side of the condition 2.d, one can use the
fact that by \eqref{Riciqpq} it mostly coincides with the one part
of the Ricci tensor of $(M,g)$.

In \cite{Boubel},  Boubel proved  theorems similar to Theorems
\ref{ThType3} and \ref{ThType4}, where conditions 2.d are changed
to  equivalent conditions on the curvature tensor. Our conditions
2.d are more precise, they can be checked using the local
coordinates and they  give the following geometric description of
the manifolds with the holonomy algebra of type 3 and 4.

Let us explain the geometric properties of the manifold with the
holonomy algebra $\g$ of type 3. Equality \eqref{lamv} shows that
the first condition of the theorem is equivalent to the equation
$\lambda=0$, where $\lambda$ is the canonically defined function
from Section \ref{secR}. In Section \ref{secR} we shown  also that
if $\lambda=0$ and a section $p$ of the parallel distribution
$\ell$ is fixed, then we obtain a sections $e$ of
$\E=\ell^\perp/\ell$ and a section $Z_\alpha$  of
$\E_\alpha\subset\E$ for each $\alpha\in K$. Condition 2.d is
equivalent to the equality
$$\pr_{\E_\alpha}e=-\frac{(n_\alpha+1)c_\alpha}{n_\alpha}J_\alpha
Z_\alpha,\quad \alpha\in K.$$ Results of Section \ref{secR} show
that if we choose $p'=\mu p$ for some non-zero function $\mu$,
then $e$ and $Z_\alpha$ change to $\mu e$ and $\mu Z_\alpha$,
respectively, i.e. the last equality does not depend on the choice
of $p$.

Manifold with the holonomy algebra $\g$ of type 4 have  the
following geometric properties. First of all, there exists a
parallel subbundle $U\subset \E_0$ of rank $n-m$. It can be
checked that if $X\in \Gamma(U)$, and a parallel light-like
vector $p$ is fixed, then the projection of $T(X)$ to
$\E=\ell^\perp/\ell$ does not depend on the choice of distribution
$E$.  We obtain $n-m$ sections $T(e_1),...,T(e_{n-m})$ of $\E$.
Condition 1 shows that these sections belong to $U^\bot$; condition
2.c. shows that the projections of these sections to $\E_\alpha$
are trivial for $\alpha\not\in K$. Thus if a vector field  $p$ is fixed,
then we obtain  sections $T(e_1),...,T(e_{n-m})$ of
$\E=\ell^\perp/\ell$ and a section $Z_\alpha$ of
$\E_\alpha\subset\E$ for each $\alpha\in K$. Condition 2.d is
equivalent to the equality
$$\pr_{\E_\alpha}T(e_s)=-\frac{(n_\alpha+1)c_{s\alpha}}{n_\alpha}J_\alpha
Z_\alpha,\quad \alpha\in K,\quad 1\leq s\leq n-m.$$ If we choose
$p'=\mu p$ for some non-zero function $\mu$, then
$T(e_1),...,T(e_{n-m})$ and $Z_\alpha$ change to $\mu
T(e_1),...,\mu T(e_{n-m})$ and $\mu Z_\alpha$, respectively, i.e.
the last equality does not depend on the choice of $p$.

Let now $(M,g)$ be an arbitrary Lorentzian manifold. Passing to
the universal covering, we may assume that $M$ is simply
connected. Results from \cite{Gal14} allow to find the
distributions defining the Wu decomposition of $(M,g)$. The
holonomy algebra of $(M,g)$ is the direct sum of the holonomy
algebras of the induced connections on the distributions
$E_\alpha$, $1\leq \alpha\leq r$ considered as  vector bundles
over $M$. If the restriction of the metric $g$ to the distribution
$E_0$ has Lorentzian signature, then the holonomy algebra of each
distribution $E_\alpha$, $1\leq \alpha\leq r$, is an irreducible
Riemannian holonomy algebra and it can be found as it is explained
at the end of Section \ref{secRiemhol}. Suppose that the
restriction of $g$ to the distribution $E_r$ is of Lorentzian
signature, then
 the holonomy algebra $\g_r$ of the induced
 connection on the distribution $E_{r}$ coincides with the holonomy
algebra of a Lorentzian manifold. If the distribution $E_{r}$
admits local non-vanishing  recurrent light-like vector fields
that are pairwise proportional in the intersection of the domains
of their definitions (or equivalently $E_{r}$ contains a parallel
subdistribution $\ell$ of isotropic lines), then the  $\g_r$ is
contained in $\simil(\rk E_r-2)$, otherwise $\g_r=\so(1,\rk
E_r-1)$. Suppose that $\g_r\subset\simil(\rk E_r-2)$. The
orthogonal part of $\g_r$ is the holonomy algebra of the induced
connection on the distribution $\ell^\bot/\ell$ (where $\ell^\bot$
is the perpendicular of $\ell$ in $E_r$), and it can be found as
in Section \ref{secorthpart}. The type of $\g_r$ can be found
using the above statements of this sections applied to local
coordinates on the integral submanifolds of the distribution
$E_r$.

\subsection{Example} Let us compare the statement of Theorem
\ref{ThType3} with the construction from~\cite{Galmetr}.

Let us fix an irreducible subalgebra $\h\subset\u(m)$, $n=2m$,
such that $\h$ contains the complex structure $J_0$ on
$\Real^{2m}$. Let $c\neq 0$ be a real number. Now we construct a
metric $g$ with the holonomy algebra
$$\g=\g^{3,\h,\varphi}=(\Real(c+J_0)\oplus\h')\zr \Real^n,\quad
\varphi(J_0)=c$$ following \cite{Galmetr}. Consider the metric
$$g=2dvdu+\sum_{i=1}^n(dx^i)^2+2Adu+H(du)^2,\quad
A=A_i(x^1,...,x^n,u)dx^i.$$ We should consider elements
$P_1,...,P_N\in \P(\h)$ such that their images generate $\h$. In
fact, it is enough to consider a single $P\in\P(\h)$: if
$\h=\u(m)$, we take $P=R^{CP^m}(\cdot,Z)$ for some non-zero vector
$Z$; if $\h\subset\u(m)$ is a symmetric Berger algebra, take
$P=R^\h(\cdot,Z)$, where $R^\h$ is the curvature tensor of a
symmetric Riemannian space with the holonomy algebra $\h$.

Let $e_1,...,e_n$ be the standard basis of $\Real^n$. Define the
numbers  $$P(e_k)e_i=P^j_{ik}e_j,\quad
a^j_{ik}=\frac{1}{3}(P^j_{ik}+P^j_{ki}),$$ then the metric is
given by \be\label{constrAH} A_j=a^j_{ik}x^ix^k,\quad
H=2vx^i\varphi(P(e_i)).\ee

Consider the conditions of Theorem \ref{ThType3}. Let $J_0e_i
=J_{0i}^je_j$. Recall that we consider the distribution $E$
spanned by the vector fields $X_i=\partial_i-A_i\partial_v$ and
the induced connection on $E$; the holonomy algebra of this
connection coincides with $\h$. We claim that the complex
structure $J$ on $E$ defined by $JX_i=J_{0i}^jX_j$ is parallel.
From \cite{Galmetr} it follows that the only nonzero Christoffel
symbols of the induced connection on $E$ are of the form
$$\Gamma^i_{uj}=P^i_{jk}x^k.$$ The condition that $\h$ commutes
with $J_0$ implies
$$P^i_{jk}J^j_{0l}=J^i_{0j}P^j_{lk}.$$
Consequently, $$\nabla_v J=\nabla_i J=0,\quad (\nabla_u
J)^i_j=\partial_u
J^i_{0j}+J^l_{0j}\Gamma^i_{ul}-J^i_{0l}\Gamma^l_{uj}=0,$$ i.e. $J$
is parallel. Next, \begin{multline*}
\nabla^jF_{lj}=\sum_{j=1}^n\partial_j(\partial_l A_j-\partial_j
A_l)=\sum_{j=1}^n\partial_j(2a^j_{lk}x^k-2a^l_{jk}x^k)=2\sum_{j=1}^n(a^j_{jl}-a^l_{jj})
\\=\frac{2}{3}\sum_{j=1}^n(P^j_{jl}+P^j_{lj}-2P^l_{jj})=-2\sum_{j=1}^nP^l_{jj}=-2(\tRic
P)_l.\end{multline*} Condition 2.d takes the form
$$\partial_i\partial_v
H=\frac{2c}{m}\sum_{l,j=1}^nP^l_{jj}J^l_{0i}.$$ Clearly, the
function
\be\label{nashliH}H=\frac{4c}{m}v\sum_{l,j=1}^nP^l_{jj}J^l_{0i}x^i\ee
satisfies this and the first conditions. From the proof of Theorem
\ref{ThType3} it follows that
$$\varphi(P(e_i))=\frac{2c}{m}\sum_{l,j=1}^nP^l_{jj}J^l_{0i}.$$
Using this and comparing \eqref{constrAH} with \eqref{nashliH}, we
see that Theorem \ref{ThType3} is in accord with the construction
from \cite{Galmetr}.

{\bf Acknowledgements.} I am grateful  to I.~M.~Anderson for
taking my attention to the problem of finding the algorithm for
computing the holonomy group of a Lorentzian manifold.  I am
thankful to D.V.~Alekseevsky  for helpful discussions.

\end{document}